# On Fuzzy Topological Spaces induced by a Given Function

Ismael Akray

Department of Mathematics \ Faculty of Science and Engineering \ Soran University

Erbil \ Kurdistan region – Iraq

ismaeelhmd@yahoo.com

*Abstract - Given a nonempty set X and a function $f: X \to X$, three fuzzy topological spaces are introduced. Some properties of these spaces and relation among them are studied and discussed.*

## 1. Introduction

Let $X$ be a nonempty set and $f: X \to X$ be a function where $f^{n+1} = f^n o f$ and $f^{-n-1} = f^{-n} o f^{-1}, \forall n \in \mathbb{N}$ ($\mathbb{N}$ is the set of natural numbers). Using this function, we introduce three fuzzy topological spaces, we study and discuss some properties of these spaces like compactness, connectedness. Finally, we give necessary and sufficient conditions under which some of these spaces coincide.

## 2. Construction of the spaces

We introduce the first fuzzy topology as follows:

Let $X$ be a nonempty set and $f: X \to X$ a function. For each $n \in \mathbb{N}$, define the fuzzy set $\mathcal{A}_n = \{(x, \mu_{\mathcal{A}_n}(x)) : x \in X\}$ where

$$\mu_{\mathcal{A}_n}(x) = \begin{cases} 1 & if\ \exists m \in \mathbb{N},\ x = f^m(x) \\ \frac{1}{n} & Otherwise \end{cases}$$

Then we have $\bigcup_{n \in \mathbb{N}} \mathcal{A}_n = X$ and $\mu_{\mathcal{A}_n \cap \mathcal{A}_m}(x) = \begin{cases} 1 & if\ \exists m \in \mathbb{N},\ x = f^m(x) \\ \frac{1}{r} & Otherwise \end{cases}$

where $r = \max\{m, n\}$. Thus, the set $\{\mathcal{A}_n : n \in \mathbb{N}\}$ is a basis for a fuzzy topology on $X$, we denote it by $\tau_1$.

**Example 2.1.** Let $f: \mathbb{N} \to \mathbb{N}$ be a function defined by $f(n) = n + 1$, then $\mathcal{A}_n = \{(x, \frac{1}{n}) : x \in \mathbb{N}\}$ and hence $\{\mathcal{A}_n : n \in \mathbb{N}\}$ is a base for the fuzzy topology $\tau_1$ on $\mathbb{N}$.

**Proposition 2.2.** In the fuzzy topological space $(X, \tau_1)$, if $\nexists m \in \mathbb{N}$, such that $x = f^m(x), \forall x \in X$, then the set $\{\mathcal{A}_n^c : n \in \mathbb{N}\}$ is a base for a fuzzy topology on $X$.

**Proof.** We have $\mu_{\mathcal{A}_n{}^c}(x) = 1 - \mu_{\mathcal{A}_n}(x) = 1 - \frac{1}{n} = \frac{n-1}{n}$, thus, $\mu_{\mathcal{A}_n{}^c \cap \mathcal{A}_m{}^c}(x) = \frac{k-1}{k}$, where $k = \min\{m,n\}$. Also, $\mu_{\cup_{n=1}^\infty \mathcal{A}_n{}^c}(x) = \sup\left\{\frac{n-1}{n} : n \in \mathbb{N}\right\} = 1$, thus $\cup_{n=1}^\infty \mathcal{A}_n{}^c = X$ and this implies that the set $\{\mathcal{A}_n{}^c : n \in \mathbb{N}\}$ is a base for a fuzzy topology on $X$.

Now, we introduce the second fuzzy topological space.

Let $X$ be a nonempty set and $f: X \to X$ be a function. Define the sets $J_0 = \cap_{n \in \mathbb{N}} f^n(X)$, $J_n = f^{n-1}(X) - f^n(X)$, $\forall n \in \mathbb{N}$. Now, for each natural number $m$, define a fuzzy subset $K_m = \left\{\left(x, \mu_{K_m}(x)\right) : x \in X\right\}$ of $X$, where $\mu_{K_m}(x) = \begin{cases} p & \text{if } x \in J_n, \ n > 0 \\ 1 & \text{Otherwise} \end{cases}$ and $p = \min\{1, \frac{n}{m}\}$. Since $\cup_{m=1}^\infty K_m = X$ and $K_m \cap K_l = K_q$ where $q = \max\{m, l\}$, the set $\mathbb{B} = \{K_m : m \in \mathbb{N}\}$ is a basis for a fuzzy topology on $X$ denoted by $\tau_2$.

**Example 2.3.** Let $f: \mathbb{N} \to \mathbb{N}$ be a function defined by $f(n) = \begin{cases} 1 & \text{if } n = 1 \\ n+2 & \text{Otherwise} \end{cases}$. Then $K_1 = \{(1,1), (2,1), (3,1), \dots\} = \mathbb{N}$, $K_2 = \{(1,1), (2, 1/2), (3, 1/2), (4,1), \dots\}$, $K_3 = \{(1,1), (2, 1/3), (3, 1/3), (4, 2/3), (5, 2/3), (6,1), (7,1), \dots\}$ and so on. Since $K_1 \supseteq K_2 \supseteq K_3 \supseteq \cdots$, and $K_1 = \mathbb{N}$, the set $\mathbb{B} = \{K_m : m \in \mathbb{N}\}$ is a base for the topology $\tau_2$ on $\mathbb{N}$.

**Lemma 2.4.** In the fuzzy topological space $(X, \tau_2)$, if $f$ is onto, then $\tau_2$ is the indiscrete fuzzy topology.

**Proof.** The proof is clear.

Finally, we introduce the third fuzzy topological space as follows:

Let $X$ be a set containing at least one element and $f: X \to X$ be a one to one function, we define a fuzzy topology on $X$ as follows:

Suppose $x_0$ is a fixed point in $X$ and $k$ is a fixed natural number, let $A(x_0) = \{y \in X : y = f^n(x_0), n \in \mathbb{N}\}$, $N_0 = \{n \in \mathbb{Z} : f^n(x_0) \in A(x_0)\}$, where $\mathbb{Z}$ is the set of integers. Then we define the following fuzzy sets for each $n$ in $N_0$ and for each $x$ in $X$:

$$C = \{(x, 1) : x \in X - A(x_0)\}$$

$$\mu_{C_n}(x) = \begin{cases} 1 & \text{if } x = f^n(x_0) \\ 1/k & \text{if } x \in A(x_0), x \neq f^n(x_0) \\ 0 & \text{Otherwise} \end{cases}$$

Now, we have $C \cap C_n = \emptyset$, $\forall n \in N_0$, $C_n \cap C_m = \left\{\left(x, \frac{1}{k}\right) : x \in A(x_0)\right\} \subseteq C_n$ and $\left(\cup_{n \in N_0} C_n\right) \cup \{C\} = X$. Then the collection $\mathbb{B} = \{C\} \cup \{C_n : n \in N_0\}$ is a base for a fuzzy topology on $X$, we denote it by $\tau_3$.

**Example 2.5.** Consider the function $f: \mathbb{N} \to \mathbb{N}$ defined by $f(n) = n + 1$ and take $x_0 = 1$, $k = 2$. Then $C = \emptyset$ and $C_n = \{(n+1, 1)\} \cup \{(m, 1/2) : m \in \mathbb{N} - \{n+1\}\}$ for all $n \in N_0 = \{0, 1, 2, \dots\}$.

**Lemma 2.6.** Let $X$ be a nonempty set endowed by the fuzzy topology $\tau_3$. Then the following are true:

1. $C = \emptyset$ iff $A(x_0) = X$.
2. The fuzzy topological space $(X, \tau_3)$ is a topological space iff $k = 1$.

**Proof.** The proof is clear.

**Theorem 2.7.** $f: (X, \tau_1) \to (X, \tau_1)$ is onto iff $f$ is open.

**Proof.** Let $f$ be onto and $x \in X$. If $x = f^m(x)$, for some $m \in \mathbb{N}$, then $\mu_{A_n}(x) = 1$ and since $f$ is onto, so $f^{-1}(x)$ is exists and one of its value is $f^{m-1}(x)$, thus $\mu_{f(A_n)}(x) = sup\{\mu_{A_n}(y): y \in X, y = f^{-1}(x)\} = \mu_{A_n}(f^{m-1}(x)) = 1$. Also, if $x \neq f^m(x)$, $\forall m \in \mathbb{N}$, then there is not exists $m \in \mathbb{N}$ such that $f^m(f^{-1}(x)) = f^{-1}(x)$. Otherwise, $f^m(x) = x$ and this is a contradiction. So, $\mu_{A_n}(f^{-1}(x)) = \frac{1}{n}$ and $\mu_{f(A_n)}(x) = sup\{\mu_{A_n}(y): y \in X, y = f^{-1}(x)\} = \frac{1}{n}$. Hence for any $x$ in $X$, we have $\mu_{A_n}(x) = \mu_{f^{-1}(A_n)}(x)$ and $f^{-1}(A_n) = A_n, \forall n \in \mathbb{N}$. Hence $f$ is open.

Now, let $f$ be open, then the image of each basic open fuzzy set is open fuzzy, hence $f(A_n)$ is open fuzzy. Suppose contrary that $f$ is not onto, so there exists $y$ in $X$ such that there is no $x$ in $X$ with $f(x) = y$ and this implies that $\mu_{f(A_n)}(y) = 0$ which means that $f(A_n)$ is not open fuzzy and this contradict with being $f$ is open.

**Theorem 2.8.** $f: (X, \tau_1) \to (X, \tau_1)$ is continuous if $f$ is one to one.

**Proof.** Let $x \in X$, if there exists a natural number $m$ such that $x = f^m(x)$, then $\mu_{A_n}(x) = 1$ and $f(x) = f^{m+1}(x)$, so $\mu_{f^{-1}(A_n)}(x) = \mu_{A_n}(f(x)) = 1$. Also, if $x \neq f^m(x), \forall m \in \mathbb{N}$, then $\mu_{A_n}(x) = \frac{1}{n}$ and since $f$ is one to one, $f(x) \neq f(f^m(x)) = f^m(f(x))$, $\mu_{A_n}(f(x)) = \frac{1}{n}$ and so $\mu_{f^{-1}(A_n)}(x) = \mu_{A_n}(f(x)) = \frac{1}{n}$. Hence $\mu_{f^{-1}(A_n)}(x) = \mu_{A_n}(f(x)), \forall x \in X$, this means that $f^{-1}(A_n) = A_n$ and $f$ is continuous.

The converse of above theorem need not be true as we see it in the following example.

**Example 2.9.** Let $f: (\mathbb{N}, \tau_1) \to (\mathbb{N}, \tau_1)$ be a function defined by $f(n) = 1$, then $f$ is not one to one and since $A_n = \{1_1, 2_{\frac{1}{n}}, 3_{\frac{1}{n}}, ...\}, \forall n \in \mathbb{N}, f^{-1}(A_n) = \{1_1, 2_1, ...\} = \mathbb{N}$ is open and $f$ is continuous.

**Theorem 2.10.** $f: (X, \tau_2) \to (X, \tau_2)$ is open iff $f$ is onto.

**Proof.** Let $f$ be onto, then by Lemma 2.4. we have $\tau_2$ is the indiscrete fuzzy topology and hence $f$ is open.

Now, suppose that $f$ is open. Since $K_1 = X, f(K_1) = f(X)$ is open. But $f(X)$ is a nonempty classical set, so $f(K_1)$ must be a classical open set and the only nonempty classical open set is $X$, hence $f(K_1) = X$ that is $f(X) = X$ and this means that $f$ is onto.

**Theorem 2.11.** $f: (X, \tau_2) \to (X, \tau_2)$ is continuous if it is onto.

**Proof.** Let $f$ be onto, then by Lemma 2.4. we have $\tau_2$ is the indiscrete fuzzy topology and hence $f$ is continuous.

The following example shows that the converse of above theorem is not true in general.

**Example 2.12.** Suppose $f: (\mathbb{N}, \tau_2) \to (\mathbb{N}, \tau_2)$ is a function defined by $f(n) = 1$, then $f$ is not onto. But, since $K_m = \{1_1, 2_{\frac{1}{m}}, 3_{\frac{1}{m}}, ...\}$ and $f^{-1}(K_m) = \{1_1, 2_1, ...\} = \mathbb{N}, \forall m \in \mathbb{N}$, so $f$ is continuous.

**Theorem 2.13.** $f: (X, \tau_3) \to (X, \tau_3)$ is open iff $f$ is onto.

**Proof.** Suppose $f$ is onto, we have to prove $f$ is open. We claim that $f(C_n) = C_{n+1}$ and $f(C) = C$. Let $x \in X$, if $x \in A(x_0)$, then there exists $r \in N_0$ such that $x = f^r(x_0)$. Thus,

$$\mu_{C_n}(x) = \begin{cases} 1 & if\ n = r \\ \frac{1}{k} & Otherwise \end{cases}, \mu_{C_{n+1}}(x) = \begin{cases} 1 & if\ n+1 = r \\ \frac{1}{k} & Otherwise \end{cases}\ and\ so$$

$\mu_{f(C_n)}(x) = \mu_{C_n}(f^{-1}(x)) = \begin{cases} 1 & \text{if } n = r-1 \\ \frac{1}{k} & \text{Otherwise} \end{cases} = \begin{cases} 1 & \text{if } n+1 = r \\ \frac{1}{k} & \text{Otherwise} \end{cases} = \mu_{C_{n+1}}(x)$. In another side, since $x \in A(x_0)$ and $f$ is onto, $f^{-1}(x_0) \in A(x_0)$, so $\mu_C(x) = 0$, $\mu_C(f^{-1}(x)) = 0$ and $\mu_{f(C)}(x) = \mu_C(f^{-1}(x)) = 0$. Thus, $\mu_{f(C)}(x) = \mu_C(x) = 0$. Hence $f(C) = C$, $f(C_n) = C_{n+1}$, $\forall x \in X$ and this implies that $f$ is open. Also, if $x \notin A(x_0)$, then $f^{-1}(x) \notin A(x_0)$ because $f$ is onto, so $\mu_C(x) = 1$, $\mu_C(f^{-1}(x)) = 1$ and $\mu_{f(C)}(x) = \mu_C(f^{-1}(x)) = 1$. Also $\mu_{C_n}(x) = 0$ and $\mu_{f(C_n)}(x) = \mu_{C_n}(f^{-1}(x)) = 0$. Thus $\mu_{C_n}(x) = \mu_{f(C_n)}(x) = 0$ and we obtain that $f$ is open.

Now, Suppose that $f$ is open, we have to prove $f$ is onto. Let $f$ be a non onto function, so there exists an element $y$ in $X$ with $y \neq f(x)$, $\forall x \in X$. If $y \in A(x_0)$, then there exists $n \in \mathbb{Z}$ such that $y = f^n(x_0)$. So $\mu_{C_n}(y) = 1$, $\mu_{C_{n+1}}(y) = \frac{1}{k}$ and $\mu_{f(C_n)}(y) = \mu_{C_n}(f^{-1}(y)) = 0$. Thus $f(C_n)$ is not open which is contradiction, since $f$ is open. But, if $y \notin A(x_0)$, then $\mu_C(y) = 1$, $\mu_{f(C)}(y) = \mu_C(f^{-1}(y)) = 0$, so $f(C) \neq C$ which means that $f(C)$ is not open and hence $f$ is not open which is contradiction. Thus $f$ is onto.

**Lemma 2.14.** In $(X, \tau_3)$, $f^{-1}(C) = C$.

**Proof.** For $x \in X$, if $x \notin A(x_0)$, then $x \neq f^n(x_0)$, $\forall n \in \mathbb{Z}$ and so $f(x) \neq f^{n+1}(x_0)$ because $f$ is one to one, so $f(x) \notin A(x_0)$ and $\mu_{f^{-1}(C)}(x) = \mu_C(f(x)) = 1$. If $x \in A(x_0)$, then $f(x) \in A(x_0)$ and $\mu_{f^{-1}(C)}(x) = \mu_C(f(x)) = 0$, hence for $x \in X$, $\mu_{f^{-1}(C)}(x) = \begin{cases} 1 & x \notin A(x_0) \\ 0 & \text{Otherwise} \end{cases} = \mu_C(x)$ and $f^{-1}(C) = C$.

**Theorem 2.15.** $f: (X, \tau_3) \rightarrow (X, \tau_3)$ is continuous if k=1.

**Proof.** Let $k = 1$, by Lemma 2.14. $f^{-1}(C) = C$, so $f^{-1}(C)$ is open. Now, we have to show that $f^{-1}(C_n)$ is open, $\forall n \in \mathbb{N}_0$., where $\mu_{C_n}(x) = \begin{cases} 1 & x \in A(x_0) \\ 0 & \text{Otherwise} \end{cases}$

Let $x \in X$, if $x \in A(x_0)$, then $\mu_{f^{-1}(C_n)}(x) = \mu_{C_n}(f(x)) = 1$. If $x \notin A(x_0)$, then $\forall n \in \mathbb{Z}$, $x \neq f^n(x_0)$ and since $f$ is one to one, $f(x) \neq f^{n+1}(x_0)$, thus $f(x) \notin A(x_0)$ and $\mu_{f^{-1}(C_n)}(x) = \mu_{C_n}(f(x)) = 0$. Therefore, $\mu_{f^{-1}(C_n)}(x) = \mu_{C_n}(x)$, $\forall x \in X$ and $f^{-1}(C_n)$ is open. Hence $f$ is continuous.

The converse of above theorem is not true, for example, consider the function $f: \mathbb{N} \rightarrow \mathbb{N}$ defined by $f(n) = 5$ and put $k = 2$ and $x_0 = 5$, then $\mu_{f^{-1}(C_0)}(n) = \mu_{C_0}(f(n)) = \mu_{C_0}(5) = 1$, so $f^{-1}(C_0) = \mathbb{N}$ and from Lemma 2.14. we have $f^{-1}(C) = C$. Hence $f$ is continuous.

### 3. Some properties of the spaces

In this section we discuss some properties for the spaces that were introduced in section two and we give necessary and sufficient conditions for the spaces to satisfy some of these properties like compactness and connectedness.

**Proposition 3.1.** $(X, \tau_1)$ is a compact space.

**Proof.** Let $\{G_\lambda : \lambda \in \Lambda\}$ be an open fuzzy cover for $X$, then $\exists \lambda_0 \in \Lambda$ such that $G_{\lambda_0} = A_1 = X$. Thus every finite open fuzzy subset of $\{G_\lambda : \lambda \in \Lambda\}$ that containing $G_{\lambda_0}$ can be consider as a finite fuzzy subcover for $X$. Hence $(X, \tau_1)$ is a compact space.

**Proposition 3.2.** $(X, \tau_1)$ is a connected space.

**Proof.** Since for every two basic open fuzzy sets $A_n, A_m$ we have either $A_n \subseteq A_m$ or $A_m \subseteq A_n$ depending on the values of $m$ and $n$ that either $m \leq n$ or $n \leq m$ respectively. Hence there are no disjoint open fuzzy sets whose union is $X$. Therefore, $(X, \tau_1)$ is connected.

**Proposition 3.3.** $(X, \tau_1)$ is not $T_0 - space$.

**Proof.** Since any two distinct fuzzy points $(x, p)$ and $(y, p)$ with equal memberships $p$ are belong to the same basic open fuzzy set $A_n$ where $p \leq \frac{1}{n}$. Hence there is no open set that contain one of the elements not the other and this implies that $(X, \tau_1)$ is not a $T_0 - space$.

**Proposition 3.4.** $(X, \tau_1)$ is regular iff there is no an element $x$ in $X$ such that $f^m(x) = x$, for some $m \in \mathbb{N}$.

**Proof.** Suppose there is no an element $x$ in $X$ such that $f^m(x) = x$, for some $m \in \mathbb{N}$, thus $\nexists x \in X$ with $\mu_{A_n}(x) = 1, \forall n \in \mathbb{N}$, so $A_n = \{(x, \frac{1}{n}): x \in X\}$ and $A_n^c = \{(x, \frac{n-1}{n}): x \in X\}$. Thus there is no element $x \in X$ and closed fuzzy set $F$ with $x \notin F$ and this means that $(X, \tau_1)$ is regular.

Now, suppose $(X, \tau_1)$ is a regular space. If there is an element $y \in X$ with $\mu_{A_n}(y) = 1$, then $y$ does not belong to any closed fuzzy set, let $F$ be one of such closed fuzzy set. But $(X, \tau_1)$ is regular, so there are disjoint open fuzzy sets one containing $y$ and the other containing $F$, but this contradict the definition of $\tau_1$ that there are no disjoint open fuzzy sets. Hence $\nexists y \in X$ with $\mu_{A_n}(y) = 1$ and this implies that $\nexists x \in X$ such that $f^m(x) = x$, for some $m \in \mathbb{N}$.

**Remark 3.5.** Since every two open fuzzy sets in $(X, \tau_1)$ have nonempty intersection, so there are no disjoint closed fuzzy sets and this implies that $(X, \tau_1)$ is a normal space.

**Proposition 3.6.** The fuzzy topological space $(X, \tau_i)$ is a Lindelof space for $(i = 1,2,3)$.

**Proof.** Since the base of the space $(X, \tau_i), i = 1,2,3$, is countable, hence every open cover has countable subcover.

**Proposition 3.7.** $(X, \tau_2)$ is a connected space.

**Proof.** Since $K_1 = X$ and $K_m \subseteq K_1, \forall m \in \mathbb{N}$. So there are no disjoint open fuzzy sets whose union be $X$. Hence $(X, \tau_2)$ is connected.

**Proposition 3.8.** $(X, \tau_2)$ is a compact space.

**Proof.** Let $A = \{G_\lambda: \lambda \in \Lambda\}$ be a open fuzzy cover for $X$. Since $K_1 = X$, there exists an element $G_{\lambda_i} \in A$ such that $K_1 \subseteq G_{\lambda_i}$ and this implies that $\{G_{\lambda_i}\}$ is a finite open subcover for $X$.

**Proposition 3.9.** $(X, \tau_2)$ is not $T_0 - space$.

**Proof.** Since each basic open fuzzy set contain every element of $X$ with nonzero membership, $(X, \tau_2)$ is not $T_0 - space$.

**Proposition 3.10.** $(X, \tau_2)$ is regular iff $f$ is onto.

**Proof.** Suppose $f$ is onto, then by Lemma 2.4. we have $\tau_2$ is the indiscrete fuzzy topology and hence $X$ is regular.

Now, let $X$ be a regular space and $f$ be a non-onto function, then $A_1 \neq \emptyset$. Take $a \in A_1$, then $\mu_{K_2}(a) = \frac{1}{2}$, so $\{(a, \frac{3}{4})\}$ is a fuzzy point in $X$ and $K_2^c$ is a closed fuzzy set not containing the fuzzy point $\{(a, \frac{3}{4})\}$, then there exist two disjoint open fuzzy sets $G$ and $H$ such that $K_2^c \subseteq G$ and $H$ containing the fuzzy point $\{(a, \frac{3}{4})\}$. But from the definition of $\tau_2$ we have no disjoint open sets, thus we have a contradiction. Hence $f$ must be onto.

**Proposition 3.11.** $(X, \tau_2)$ is a normal space.

**Proof.** Since the basic open fuzzy sets have the property that $X = K_1 \supseteq K_2 \supseteq \cdots$, so the closed fuzzy sets have the property that $\emptyset = K_1^c \subseteq K_2^c \subseteq \cdots$, hence there are no disjoint closed sets and consequently $X$ is normal.

**Proposition 3.12.** $(X, \tau_3)$ is connected iff $C = \emptyset$.

**Proof.** Suppose $X$ is connected, we have to prove $C = \emptyset$. If $C \neq \emptyset$, then take $A = \bigcup_{n \in N_0} C_n$ and hence $C$ and $A$ are two disjoint nonempty open fuzzy sets whose union is $X$ and this is a contradiction. Thus, $C = \emptyset$.

Now, let $C = \emptyset$, since for any basic open fuzzy sets $K_m$ and $K_n$ we have $K_m \cap K_n = \{(x, \frac{1}{k}): x \in X\}$ which is nonempty. Hence $X$ is connected.

**Proposition 3.13.** $(X, \tau_3)$ is a compact iff $N_0$ is finite.

**Proof.** Suppose $N_0$ is a finite set say $N_0 = \{n_1, n_2, \dots, n_m\}$ and let $\{G_\lambda: \lambda \in \Lambda\}$ be an open fuzzy cover for $X$, then there exist $G_{\lambda_1}, G_{\lambda_2}, \dots, G_{\lambda_m}, G_{\lambda_0}$ in $\{G_\lambda: \lambda \in \Lambda\}$ such that $C_{n_i} \subseteq G_{\lambda_i}$ for $i = 1, 2, \dots, m$ and $C \subseteq G_{\lambda_0}$. Hence $\{G_{\lambda_0}, G_{\lambda_1}, G_{\lambda_2}, \dots, G_{\lambda_m}\}$ is a finite open fuzzy subcover for $X$ and this means that $X$ is compact.

Now, let $X$ be a compact space and $N_0$ infinite set, then $\{C_n: n \in N_0\} \cup \{C\}$ is an open fuzzy cover for $X$ that have no finite subcover, so $X$ is not compact and this contradict with our assumpsion. Thus $N_0$ must be finite.

**Proposition 3.14.** $(X, \tau_3)$ is $T_0 - space$ iff $C$ has atmost one element and $f(x_0) = x_0$.

**Proof.** Suppose $C$ has atmost one element and $f(x_0) = x_0$. If $C = \emptyset$, then $X = \{x_0\}$ and $X$ is $T_0 - space$. If $C \neq \emptyset$, then there exists $y \neq x_0$ in $X$ such that $C = \{y\}$, that is $X = \{x_0, y\}$ and there is an open fuzzy set $C$ contain $y$ but not $x_0$. Hence $X$ is $T_0 - space$.

Now, let $(X, \tau_3)$ be a $T_0 - space$. If $f(x_0) \neq x_0$, then there exists an element $y$ in $X$ such that $f(x_0) = y$ and this means that $x_0 \neq y$ and there is no open fuzzy set contains only one of them, so $X$ is not $T_0 - space$ which is contradiction. Also, if $C$ contains two elements, then it makes $X$ to be non $T_0 - space$. Hence must $f(x_0) = x_0$ and $C$ contains atmost one element.

**Proposition 3.15.** $(X, \tau_3)$ is regular iff $f(x_0) = x_0$.

**Proof.** Suppose $f(x_0) = x_0$ and $x \in X$. Let $F$ be a nonempty closed fuzzy set not containing $x$. We have two cases; first, if $= x_0$, then $F = C_0^c$, so there exist two disjoint open fuzzy sets $C$ and $C_0$ such that $F \subseteq C$ and $x \in C_0$. Second, if $\neq x_0$, then $x \in C$ and $F = C^c$, so there exist two disjoint open fuzzy sets $C$ and $C_0$ such that $F \subseteq C_0$ and $x \in C$. In both cases we conclude that $X$ is regular.

Now, let $X$ be a regular space and $f(x_0) \neq x_0$. Put $x = x_0$ and $y = f(x_0)$, then for the closed fuzzy set $C_0^c$ not containing $x$, there exist two disjoint open fuzzy sets $G$ and $H$ such that $x \in G$ and $C_0^c \subseteq H$. Since $x \notin H$, $H = C$, but $\mu_{C_0^c}(y) = \frac{k-1}{k}$ and $\mu_C(y) = 0$, so $C_0^c \nsubseteq C$ which is a contradiction. Hence $x = y$, that is $f(x_0) = x_0$.

**Proposition 3.16.** $(X, \tau_3)$ is a normal space.

**Proof.** To prove $X$ is normal, we have two cases;

1. If $f(x_0) = x_0$ and $C \neq \emptyset$, then there are only two nonempty disjoint closed fuzzy sets $C^c$ and $C_0^c$, so there are two disjoint open fuzzy sets $C$ and $C_0$ such that $C_0^c \subseteq C$, $C^c \subseteq C_0$ and hence $X$ is normal.
2. If either $f(x_0) \neq x_0$ or $C = \emptyset$, then we have no disjoint nonempty closed fuzzy sets, thus $X$ is normal.

## 4. Relation among the spaces

In this section we study the necessary conditions for some of the spaces to be coincide. For this purpose we have the following theorems.

**Theorem 4.1.** The two fuzzy topologies $\tau_1$ and $\tau_2$ are equal iff for every $x \in X$, there exists $m \in \mathbb{N}$ such that $f^m(x) = x$.

**Proof.** Assume that for every $x \in X$, there exists $m \in \mathbb{N}$ such that $f^m(x) = x$, then $K_m = X$, $\forall m \in \mathbb{N}$ and $\tau_2$ is the indiscrete fuzzy topology. Furthermore, we have from the assumption that for every $x \in X$, $\mu_{A_n}(x) = 1$, so $\tau_1$ is the indiscrete fuzzy topology. Therefore, $\tau_1 = \tau_2$.

Now, suppose that $\tau_1 = \tau_2$ and according to the definitions of $\tau_1$ and $\tau_2$, we have for every $x \in X$, $\mu_{A_n}(x) = \begin{cases} 1 & if\ \exists m \in \mathbb{N},\ x = f^m(x) \\ \frac{1}{n} & Otherwise \end{cases}$ and $\mu_{K_m}(x) = \begin{cases} min\{1, \frac{i}{m}\} & if\ x \in J_i,\ i > 0 \\ 1 & Otherwise \end{cases}$. By contrary that if there exists $x$ in $X$ such that for every natural number $m$, $f^m(x) \neq x$, then $\mu_{A_n}(x) = \frac{1}{n}$, $\mu_{K_m}(x) = min\{1, \frac{i}{m}\}$ which they are not equal. Therefore, for every $x \in X$, there exists $m \in \mathbb{N}$ such that $f^m(x) = x$ and this completes the proof.

**Theorem 4.2.** The two fuzzy topologies $\tau_1$ and $\tau_3$ are never be equal.

**Proof.** From the fuzzy topology $\tau_3$, $\mu_C(x_0) = 0$, but there is no open fuzzy set $G$ in $\tau_1$ and $x$ in $X$ such that $\mu_G(x) = 0$. Thus $\tau_1 \neq \tau_3$.

**Theorem 4.3.** The two fuzzy topologies $\tau_2$ and $\tau_3$ are equal iff $f$ is onto, $A(x_0) = X$ and $k = 1$.

**Proof.** Suppose $f$ is onto, $A(x_0) = X$ and $k = 1$, then $A_n = \emptyset$, $\forall n > 0$ and $A_0 = X$, so $\mu_{K_m}(x) = 1$, $\forall x \in X$, $\forall m \in \mathbb{N}$ and hence $\tau_2$ is the indiscrete fuzzy topology. Also, since $A(x_0) = X$, $C = \emptyset$ and $\mu_{C_n}(x) = 1$. Thus $C_n = X$, $\forall n \in N_0$ and $\tau_3$ is the indiscrete fuzzy topology. Hence $\tau_2 = \tau_3$.

Now, let $\tau_2 = \tau_3$, since for each $x$ in $X$, we have $\mu_{K_m}(x) = \begin{cases} min\{1, \frac{i}{m}\} & if\ x \in J_i,\ i > 0 \\ 1 & Otherwise \end{cases}$ and $\mu_{C_n}(x) = \begin{cases} 1 & if\ x = f^n(x_0) \\ 1/k & if\ x \in A(x_0), x \neq f^n(x_0) \\ 0 & Otherwise \end{cases}$

So they are equal if $\frac{i}{m} \geq 1$, $C = \emptyset$ and $\frac{1}{k} = 1$ that is $i \geq m$, $A(x_0) = X$ and $k = 1$. But for $K_2$, we have $i \geq 2$ and this implies that $J_1 = \emptyset$. Hence $f$ is onto, $A(x_0) = X$ and $k = 1$.